\documentclass[a4paper, 10pt, conference]{ieeeconf}      % Use this line for a4 paper

\IEEEoverridecommandlockouts                              % This command is only needed if 
\newcommand{\ncom}{\newcommand}
\ncom{\beqn}{\begin{eqnarray*}} \ncom{\eeqn}{\end{eqnarray*}}
\ncom{\beq}{\begin{eqnarray}} \ncom{\eeq}{\end{eqnarray}}

\ncom{\doq}{\nabla _q} \ncom{\dop}{\nabla _p} \ncom{\dopt}{\nabla
_{\tilde{p}}}

\usepackage{amsmath}
\usepackage{amsfonts}
\usepackage{amssymb}
\usepackage{graphicx}
\newtheorem{theorem}{Theorem}
\newtheorem{remark}[theorem]{Remark}
\newtheorem{proposition}[theorem]{Proposition}
\newtheorem{example}[theorem]{Example}

\newtheorem{defn}{Definition}

\title{\LARGE \bf
On Power Balancing and Stabilization for a Class of infinite-dimensional systems}

\author{Ramkrishna Pasumarthy, Krishna Chaitanya Kosaraju, Addarsh Chandrasekar%
\thanks{ Ramkrishna Pasumarthy and Krishna Chaitanya Kosaraju are with Department of Electrical Engineering, IIT Madras, India.
        {\tt\small ramkrishna@ee.iitm.ac.in, ee1d055@ee.iitm.ac.in}.
        Addarsh Chandrasekhar is an undergrad student at the Department of Mechanical Engineering, IIT Madras, India. 
}
}
\begin{document}

\maketitle
\thispagestyle{empty}
\pagestyle{empty}

%%%%%%%%%%%%%%%%%%%%%%%%%%%%%%%%%%%%%%%%%%%%%%%%%%%%%%%%%%%%%%%%%%%%%%%%%%%%%%%%
\begin{abstract}
In this paper we present control of infinite-dimensional systems by power shaping methods, which have been used extensively for control of finite dimensional systems. Towards achieving the results we work within the Brayton Moser framework, by using the system of transmission line as an example and derive passivity of the system with respect to the boundary voltages and derivatives of current at the boundary. We then solve the stabilization problem by interconnecting the system through a finite-dimensional controller and generating Casimirs for the closed-loop system. Finally we explore possibility of generating other alternate passive maps. 
\end{abstract}

%%%%%%%%%%%%%%%%%%%%%%%%%%%%%%%%%%%%%%%%%%%%%%%%%%%%%%%%%%%%%%%%%%%%%%%%%%%%%%%%
\section{INTRODUCTION}
The concept of passivity based control has been well incorporated towards solving stabilization problems arising in nonlinear systems. In particular the Hamiltonian and Lagrangian frameworks naturally incorporate the concept of passivity with the total energy serving as the storage function and port variables which are power conjugate (for example voltages and currents, forces and velocities). This has resulted in a control technique usually referred to as control by ``energy shaping" of physical systems \cite{l2gain,geoplex,ortaut}. In some cases the ``natural" power conjugate port variables do not necessarily help in achieving the control objectives due to the ``dissipation obstacle", motivating the search for alternate passive maps or other possible storage functions \cite{guido,eloisa,ortegapower}. One possibility is to use the so called mixed-potential function as the storage function, which has led to a control method called the power shaping technique. This mixed potential function arises from the Brayton Moser framework for modeling of RLC circuits \cite{BM01,BM02}. Several of these energy and power shaping techniques have been reported in literature for systems which are lumped in nature.
Towards control of distributed parameter systems by energy shaping, few results have been reported in \cite{hugo,geoplex}. However, to the best of our knowledge, the extension of control of infinite-dimensional systems via power shaping methods hasn't been addressed so far. Results on stability in the case of transmission line as a subsystem have been reported in \cite{BM03}.

In this paper we use the results on stability of infinite dimensional systems reported in \cite{BM03} and study problems on control of infinite dimensional systems via power shaping methods. The use of mixed potential functions as the storage functions results in passivity properties with respect to the boundary voltages and derivatives of current at the boundary. The controller is a finite dimensional system interconnected to the boundary of the transmission line and control strategy and stabilization will rely on finding Casimir functions of the closed-loop systems. Similar results on boundary control of infinite-dimensional systems, via energy shaping methods, have  been reported in \cite{hugo,lzm}.

The rest of the paper is organised as follows: In Section II we briefly recall the Brayton Moser (BM)  framework for RLC circuits and present an example showing control by power shaping with example of an RLC circuit. In Section III we present the example of the transmission line modelled in the BM framework and generate admissible pairs to infer stability of the transmission line under the conditions of zero boundary energy flow. We also show how to incorporate boundary conditions in the BM framework for the transmission line. These results are motivated by those presented in \cite{BM03,jeltcsm}. In Section IV we present a control problem of stabilizing the transmission line at a non-zero equilibrium. This is achieved by interconnecting the transmission line to a controller via the boundary. To achieve the desired control objective we rely on shaping the power of the system by generating Casimir functions for the closed-loop system. In Section V we present the possibility of generating further alternate passive maps, where the storage function now has the dimensions of power/time and passivity is inferred via the port variables which are now the derivatives of the currents and the derivatives of the voltages at the boundary. Finally we end with some concluding remarks.

\section{Brayton Moser Framework}
Brayton and Moser in the early sixties \cite{BM01, BM02} showed that the dynamics of a
class (topologically complete) of nonlinear $RLC$-circuits can be
written as
\begin{align}
A (i_L,v_c)
 \begin{bmatrix}
\frac{di_L}{dt} \\ \frac{dv_c}{dt}
\end{bmatrix} =
\begin{bmatrix}
\frac{\partial P}{\partial i_L} \\
\frac{\partial P}{\partial v_C}
\end{bmatrix} +
\begin{bmatrix}
B_{E_c}^T E_c \\
-B_{J_c}^T J_C
\end{bmatrix}
\label{eq_bm_fd1}
\end{align}
where $A(i_L,v_C) = \text{diag} \{L(i_L), -C(v_C) \}$ and $i_L$ the vector
of currents through inductors, $v_C$ vector of capacitor voltages,
$L(i_L)$ the inductance matrix, $C(v_C)$ the capacitance matrix,
$B_{E_c}, B_{J_c}$ the matrices containing the elements $\{-1, 0 , 1
\}$ decided by Kirchoff's voltage and current laws. $E_C, J_C$ are
respectively the controlled voltage and current sources. $P$ is
called the mixed potential function defined by
\[
P(i_L,v_C) = F(i_L) - G(v_C) + i_L^T \gamma\; v_c
\]
where $x = (i_L, v_C)$ the system states. Here $F$ is the content of
all the current controlled resistors, $G$ is the co-content of all
voltage controlled resistors. Matrix $\gamma$ contains elements
 $\{-1, 0 , 1 \}$ depending on the network topology.
Computing the time derivative of $P$ along the trajectories of
(\ref{eq_bm_fd1}) we have
\[
\dot P = \dot x ^T (A(x) + A^T(x)) \dot x + u^T y,
\]
where, $u=\left(E_c,\;J_c\right)$ and
$y=\left(-B_{E_c}\frac{di_L}{dt},\;B_{J_c}\frac{dv_C}{dt}\right)$.\\
From the above expression we can conclude that the system is
passive if $(A(x) + A^T(x))\le 0$.

In case $(A(x) + A^T(x))\le 0$ is not satisfied, then it is possible
to find a new $(\tilde A, \tilde P)$ called an ``admissible pair", (refer \cite{guido,ortegapower})
satisfying $(\tilde A(x) + \tilde A^T(x))\le 0$.  The dynamics can
equivalently be written as
\begin{align}
 \tilde A \begin{bmatrix}
\frac{di_L}{dt} \\ \frac{dv_c}{dt}
\end{bmatrix} =
\begin{bmatrix}
\frac{\partial \tilde P}{\partial i_L} \\
\frac{\partial \tilde P}{\partial v_C}
\end{bmatrix} +
\begin{bmatrix}
B_{E_c}^T E_c \\
-B_{J_c}^T J_C
\end{bmatrix}
\label{eq_bm_fd_adm1}
\end{align}
\begin{remark}
Contrast to the case where the total energy of the system serves as
the storage function and derive passivity with respect to input -
output variables which are power conjugate, for example the voltage
and currents \cite{l2gain}. In the case of using the mixed potential as the
storage function we derive passivity with respect to controlled
voltages, and the derivatives of currents or the controlled currents
and the derivatives of the voltages.
\end{remark}

\subsection{Stabilization by Power Shaping}
We consider the problem of stabilizing a given system at any desired
equilibrium point. What if $\tilde P$ does not have a minimum
at the desired equilibrium? The idea is then to generate dynamical
invariants independent of the Hamiltonian called Casimir functions
$\mathcal C : X \rightarrow \mathbb R$ such that $\frac{d \mathcal
C}{dt} = 0$ along the system trajectories. Here $X \subset \mathbb R^n$ denotes the state space of the system. We then construct a
Lyapunov function as the sum of the power function $\tilde P$ plus
the Casimir function. For a comprehensive exposure to various energy
and power shaping techniques we refer to \cite{guido,ramijc,eloisa,ortegapower}.

Power shaping
stabilization is a method where the storage function, instead of the
total stored energy of the system, is a function related with the
power of the circuit, called a mixed potential function. In the
context of electrical networks, the passivity property is now
established with respect to voltage and derivative of current, or
current and derivative of voltage. This modeling approach was first
proposed by Brayton-Moser \cite{BM01} for non linear RLC circuits.
In this framework the control objective is now achieved via
assigning a new power function to the system and hence is referred
to as stabilization by``power shaping"
\cite{eloisa,geoplex,ortegapower}. This method has natural
advantages over practical drawbacks of energy shaping methods like
speeding up the transient response (as derivatives of currents and
voltages are used as outputs) and also help overcome the
``dissipation obstacle" \cite{geoplex}. We now present results on control by power shaping with help of an example of a simple
RLC circuit presented in \cite{eloisa}.
\begin{figure}
[ptbh]
\begin{center}
\includegraphics[height=1.1 in,width=2.5in]%
{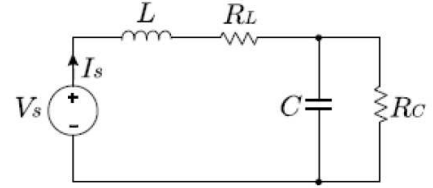}%
\caption{A RLC circuit}
\label{fig:rect}%
\end{center}
\end{figure}
\begin{example}
Consider the circuit of Figure \ref{fig:rect}. In the Brayton-Moser
framework, the dynamics of the circuit can be written as
\begin{align*}
A\begin{bmatrix}
\frac{d}{dt}i_L \\
\frac{d}{dt} v_C
\end{bmatrix} =
\begin{bmatrix}
\frac{\partial P}{\partial i_L} \\
\frac{\partial P}{\partial v_C}
\end{bmatrix} +
\begin{bmatrix}
-1 \\
0
\end{bmatrix}
v_S
\end{align*}
with $A = \text{ diag } \{-L, C \}$, $i_{L}$ is the current through
the inductor and $v_{C}$ the voltage across the Capacitor. $P$ is
the mixed potential function (power function) given by
\[
P =- \frac{R_C}{2} \left ( \frac{v_C}{R_C} - i_L\right
)^2+\frac{1}{2} (R_L + R_C)i_L^2
\]
The equilibrium points $(i_L^{\ast}, v_C^{\ast})$  of
the system are given as
\[
v_C^{\ast} = R_C i_L^{\ast}, ~~ i_L^{\ast} = \frac{v_S^{\ast}}{R_C +
R_L}
\]
To deduce the desired passivity property (with respect to $v_S$ and
$\frac{di_S}{dt}$), we need to find admissible pairs $\tilde A,
\tilde P$ such that
\begin{align*}
\tilde A\begin{bmatrix}
\frac{d}{dt}i_L \\
\frac{d}{dt} v_C
\end{bmatrix}
=
\begin{bmatrix}
\frac{\partial \tilde P}{\partial i_L} \\
\frac{\partial \tilde P}{\partial v_C}
\end{bmatrix} +
\begin{bmatrix}
-1 \\
0
\end{bmatrix}
v_S
\end{align*}
As shown in \cite{eloisa}, the following choice of the admissible
pair
\begin{align*}
\tilde A & = \begin{bmatrix}
-L & 2R_CC \\
0 & -C
\end{bmatrix} \\
\tilde P & = \frac{R_C}{2} \left ( \frac{v_C}{R_C} -i_L\right )^2 +
\frac{1}{2} (R_L + R_C)i_L^2,
\end{align*}
together with the condition $L\geq R_C^2 C$, results in the desired
dissipation inequality as $ \dot {\tilde {P}} \le \frac{d
i_S}{dt}v_S. $ Further, we can achieve the required stabilization
via the controlled voltage \cite{eloisa,geoplex}
\[
v_S = -K(i_L - i_L^{\ast}) + (R_L + R_C)i_L^{\ast}
\]
with $K\ge 0$ a tuning parameter. This control law globally
stabilizes the system with Lyapunov function
\begin{align}
\tilde P_d = \frac{R_C}{2} \left ( \frac{v_C}{R_C} -i_L\right )^2 +
\frac{1}{2} (R_L + R_C + K)(i_L^2 - i_L^{\ast})^2 \label{eq_newp}
\end{align}
Thus we have assigned a new power function to the controlled system
which achieves the desired stability property.
\end{example}

\section{The transmission line}
In this section, we present an example of an infinite-dimensional system, the transmission line, in the Brayton Moser (BM) framework, which originally appeared \cite{BM03}. Based on the arguments in the previous section, we aim to generate admissible pairs for the system of equations governing the transmission line with time varying boundary conditions, resulting in a non zero boundary energy flow, and also explore alternate passivity properties for the transmission line in the BM framework, which will be used for control purposes. The results presented in this section are motivated by results that appeared in \cite{BM03, ja}.

The dynamics of a transmission line are described by the following telegraphers equations
\begin{align}
L\frac{d}{dt}i(z,t) = -\frac{\partial}{\partial z}v(z,t) -Ri(z,t) \nonumber \\
-C \frac{d}{dt}v(z,t) = \frac{\partial}{\partial z}i(z,t) + G v(z,t),
\end{align}
where $0<z<1$ represents the spatial domain of the transmission line with $L,C$ respectively denoting the specific inductance and capacitance, which in this case is assumed to be independent of the spatial variable $z$. Further, $R,G$  denote the specific resistance and conductance respectively, which are also assumed to be independent of $z$. Finally, we denote the boundary voltages and currents as $v_0(t) = v(0,t), i_0(t) = i(0,t)$ and $v_1(t) = v(1,t), i_1(t) = i(1,t)$

The dynamics can be re-written in the following BM form:
\begin{align}
\begin{bmatrix}
L& 0 \\
0 & -C
\end{bmatrix}
\begin{bmatrix}
\frac{\partial}{\partial t}i(z,t)\\
\frac{\partial }{\partial t}v(z,t)
\end{bmatrix}
=
\begin{bmatrix}
\delta_i  P \\
\delta_v  P
\end{bmatrix}
\label{eq_tlbm1},
\end{align}
where $P$ is  a mixed potential functional \cite{BM01,BM02} given by
\begin{align}
P & = \int_0^1 \left (\left( -\frac {\partial v}{\partial z}\right ) i + \frac{1}{2} G v^2 - \frac{1}{2}R i^2 \right )dz \nonumber \\
& = \int_0^1 \bar P \left (i,v, \frac {\partial v}{\partial z} \right) dz 
\label{eq_ptl_1}
\end{align}
and $\delta$ denotes the variational derivative. The above expression can be treated as the infinite-dimensional version of (\ref{eq_bm_fd1}). Boundary conditions can be incorporated as follows: The spatial domain is of finite length $0 < z < 1$. Suppose at $z=0$ is a voltage source $E$ in series with a resistor $R_0$ and at $z=1$ is connected a capacitor of value $C_1$ with a resistance $R_1$ in parallel. This would result in the following dynamics at the boundary:
\begin{align}
E & = v_0 + i_0 R_0,& z=0 \nonumber \\
i_1 & = C_1\frac{dv_1}{dt} + \frac{v_1}{R_1}, & z=1
\label{eq_tlbc1}
\end{align}
The new mixed potential incorporating the boundary conditions can now be written as
\begin{align}
\mathcal P = P +  P^0 + P^1
\label{eq_P}
\end{align}
where $P^0 = (E-v_0 )i_0 - \frac{1}{2}i_0 ^2R_0$ and $P^1 =
\frac{1}{2}\frac{ v_1^2}{R_1}$. Now, in addition to equations
(\ref{eq_tlbm1}) we have the following conditions governing the dynamics at the boundary:
\begin{align}
&0 =\delta_{i_0} \mathcal P =  (E-v_0)-i_0 R_0 \nonumber \\
&\delta_{v_1} \mathcal P = \frac{v_1}{R_1} = -i_1  -C_1
\frac{dv_1}{dt} \label{eq_b1}
\end{align}
The transmission line, which is an infinite-dimensional system and whose dynamics are governed by equations (\ref{eq_tlbm1}) together with the boundary dynamics given by (\ref{eq_tlbc1}) can be written in a unified way as
\begin{align}
A \frac{\partial u}{\partial t} = J \frac{\partial u}{\partial z} + G
\end{align}
where
\begin{align*}
& u = \begin{bmatrix}
i(z,t) &
v(z,t) &
i(0,t) &
v(0,t) &
i(1,t) &
v(1,t)
\end{bmatrix}^T \\
& A = \text{diag} \{A_0, A_1, A_2 \}, J = \text{diag}\{ J_0 , 0 , 0\}\\
& G = [G_0~ G_1 ~G_2]^T
\end{align*}
and
\begin{align}
&A_0 = \text{diag} \{L, -C \}, A_1 = 0, A_2 = \text{diag}\{0 , -C_1 \} \nonumber \\
&G_0 = [-Ri~ Gv]^T, G_1 = [E-V_0-R_0i_0~ 0]^T, \nonumber & \\ &G_2 = [0 ~ -i_1 + \frac{v_1}{R_1}]^T 
J_0 =
\begin{bmatrix}
0 & -1 \\
1 & 0
\end{bmatrix}
\end{align}
The system, can more generally be written in canonical form as
\begin{align}
A(u) \frac{\partial u}{\partial t} =\delta_u\mathcal P
%\frac{\partial \mathcal P}{\partial u}
\label{eq_can}
\end{align}
where the functional $\mathcal P$ can be
written as
\[
\mathcal P(u) = \int_0^1 \bar P\left (u,\frac{\partial}{\partial
z}u\right ) dz + \phi(u)\mid_{z=0} + \left.\psi(u)\right|_{z=1}
\]
 In our case
 $\delta_u\mathcal P$ is computed as
  \begin{align}
\delta_u\mathcal P =
 \begin{bmatrix}
 -\frac {\partial v}{\partial z}-Ri \\
 Gv+\frac {\partial i}{\partial z}\\
 E-v_0-R_0i_0\\
 0\\
 0\\
 \frac{v_1}{R_1}-i_1
 \end{bmatrix}
 \end{align}
 with $\mathcal P$ as in (\ref{eq_P}).
 The systems of the form (\ref{eq_can}) are usually referred to as mixed finite and infinite-dimensional systems and also sometimes called mixed initial boundary value problems.
\subsection{Stability}
\begin{defn}
$u\in C_1[0,1]$ means the first $n$ components of $u$ are continuous
in $[0,1]$ and have a piecewise continuous derivative.
\end{defn}
\begin{defn}
By Hilbert space $H_1$ we mean the completion of the space of
functions $u \in C_1[0,1]$ with the norm  $||u||_1$ given by
\beqn
||u||_1^2&=&(u(0),u(0))+(u(1),u(1))+\int_0^1(u,u)dz\\
&& +\int_0^1\left(\frac{\partial u}{\partial z},\frac{\partial
u}{\partial z}\right)dz
\eeqn
where $(\cdot,\cdot)$ represents the usual inner product.
\end{defn}
\begin{defn}
By Hilbert space $H_0$ we mean the completion of the space of
functions $u \in C_1(0,1)$ with the norm  $||u||_0$ given by
\end{defn}
\beqn
||u||_0^2&=&(u(0),u(0))+(u(1),u(1))+\int_0^1(u,u)dz
\eeqn
To infer stability of such systems, which are mixed finite and infinite-dimensional in nature, we refer to results presented in \cite{BM03} which we briefly recall here. For comprehensive exposure to the results we refer the readers to \cite{BM03} and references therein.
\begin{defn}
Denote by $F_t,~~ t \in (-\infty, \infty)$ a one parameter family of operators which define a dynamical system in some metric space $R$. Let $L(p)$ denote the limit set of a trajectory $F_t(p), ~ p\in R$ and $M$ be an invariant set of the dynamical system. If for each $p\in R, ~ F_t(p), ~ t \ge0$ is compact in $R$, then every limit set $L(p)$ is non empty. Moreover, if there exists an invariant set $M$ containing $L(p)$ for all $p\in R$, then $M$ is said to be completely (uniformly asymptotically) stable.
\end{defn}
\begin{defn}
The equilibrium set $\mathcal E$ of (\ref{eq_can}) is the set of
functions $u \in H_1$ such that $\left.|| \delta_u\mathcal
P||\right._0 =0$
\end{defn}
Towards proving stability we invoke the following theorem from \cite{BM03}:
\begin{theorem}
Let
\[
 A\frac{\partial}{\partial t}u = \delta_u\mathcal P = J \frac{\partial}{\partial z} u+ G, ~~ z \in (0,1), ~ t>0
 \]
 Here $A, J, G$ depend continuously on $z$ and $u$. If the following conditions are satisfied then the equilibrium set $\mathcal E$ is completely stable
 \begin{enumerate}
 \item $\delta_u P \in H_0$ and is continuous for $u \in H_1$ a Hilbert space
 \item $J_0$ is invertible
 \item There exists a functional $\tilde P (u)$ called the Lyapuniv functional, which is a continuous functional for $u \in H_0$ another Hilbert space, and is such that
 \begin{enumerate}
 \item $\tilde {\mathcal P}(u) \ge -b > -\infty, ~ b$ some real number
 \item $\tilde {\mathcal  P}(u) \rightarrow \infty \Leftrightarrow ||u||_1 \rightarrow \infty $
 \item $\dot{\tilde {\mathcal P}} \le -K ||\delta_u \mathcal  P ||^2_0 \le 0$, $K$ some positive number.
 \end{enumerate}
 \end{enumerate}
 \label{thm}
\end{theorem}

\subsection{Stability of the Transmission line}
We present several cases for testing the stability of the transmission line. To begin with let us assume that there is no net power flow through the boundary of the system.
Consider now that function $ P $ as in equation (\ref{eq_ptl_1}), which has dimensions of power. The rate of this
function is computed as
\begin{align*}
& \dot { P}  =\int_0^1\left (  L \left (\frac{di(z,t)}{dt} \right)^2 - C \left ( \frac{dv(z,t)}{dt}\right )^2 \right )dz
\end{align*}
It can be easily seen that the right hand side of the above equation is not sign definite, and hence $ P$ does not serve as a Lyapunov functional to infer any kind of stability (or for that matter passivity) properties of the system.  We thus need to look for other possible Lyapunov functionals $\tilde P$, or in other words ``admissible pairs" $\tilde{\mathcal A}, \tilde{\mathcal P}$ as in the case of finite dimensional systems (refer equation (\ref{eq_bm_fd_adm1})) which can prove stability of the system.
%The equations of the transmission line can be written in a more generally as a pseudo-gradient system of the form
%\[
%M \frac{\partial u}{\partial t} = \delta_u P
%\]
%with respect to the indefinite (pseudo-Riemannian) metric which in the case of the transmission line takes the form
%\[
%M =
%\begin{bmatrix}
%-L & 0 \\
%0 & C
%\end{bmatrix}
%\]
The admissible pairs $\tilde {\mathcal P}$ and $\tilde{\mathcal A}$ should be  such that in the new pair we still retain the dynamics of (\ref{eq_tlbm1}) in case of the transmission line or equation (\ref{eq_can}) for the more general case, of course with zero boundary energy flow.  Moreover the admissible pair should be such that the symmetric part of $\tilde A$ is negative semi-definite. This can be achieved in the following way. Let
\[
{\tilde P}(u) = \lambda { P} + \frac{1}{2}\int_0^1 \delta_u  P(u) \cdot K \delta_u P(u))dz
\]
with $\lambda$ being an arbitrary constant and $K$ a symmetric matrix.
The aim is to find values of $\lambda$ and $K$ such that
\[
\dot {\tilde P} = \int_0^1 \frac{\partial u}{\partial t} \cdot \tilde{\mathcal A} \frac{\partial u}{\partial t} ~ dz \le 0.
\]
A possible choice, of course not unique, is choosing $\lambda = -1$ and $K = \text{diag } (0, 2/G)$ yielding a new mixed potential function of the form
\begin{align}
{\tilde P} & = \int_0^1 \left  (\frac{1}{2} R i^2 + \frac{1}{2} Gv^2 + v \frac{\partial i}{\partial z}  + \frac{1}{G} \frac{\partial i}{\partial z}^2 \right )dz  \nonumber \\
& = \int_0^1 \frac{1}{2G}\left ( \frac{\partial i}{\partial z} + G v\right )^2dz + \int_0^1 \frac{1}{2G}\left ( \frac{\partial i}{\partial z}\right )^2dz \nonumber \\
& + \int_0^1 \frac{1}{2}Ri^2 dz
\label{eq_newP}
\end{align}
This new mixed potential function $\tilde P$ satisfies condition 3b of Theorem \ref{thm}.
The  $\tilde A$ is then verified to be
\begin{align}
\tilde A =
\begin{bmatrix}
-L & \frac{2C}{G} \frac{\partial}{\partial z}\\
0 & -C
\end{bmatrix}
\label{eq_newA}
\end{align}
It is easy to check that $i(z,t) = 0 = v(z,t)$ is an equilibrium point of the system.
If the symmetric part of $\tilde A$ is negative semi definite, then the new mixed potential (\ref{eq_newP}) will be non increasing along the trajectories of the system, provided $\sqrt{C/L}/G ||\frac{\partial}{\partial z}|| \le 1$, thus satisfying condition $3c$ of Theorem \ref{thm} and hence inferring complete stability of the system.
%\subsection{stability of the systems under boundary conditions}
%Consider the transmission line system described by equations () together with () as the boundary conditions. Consider the new mixed potential function

\subsection{Alternate passive maps for the transmission line}
In this section we derive passivity properties of the transmission line for non zero boundary conditions. In the previous section on finite dimensional systems, we have shown how with the help of admissible pairs, we can generate a new mixed potential function (\ref{eq_bm_fd_adm1}) and infer passivity with respect to port variables which are the voltages and derivatives of current. In the classical energy based modeling of the transmission line \cite{stokes,hugo} where the total energy of the system is used as the storage functional, we infer passivity with respect to the boundary voltages and currents. Here we show how we can make use of the  power function (\ref{eq_newP}) to derive passivity properties with port variables as the boundary voltages and derivatives of current (or vice versa)
Consider the formulation of the transmission line as in equations (\ref{eq_tlbm1}) together  with boundary conditions as in equations (\ref{eq_b1}). For simplicity we modify the boundary conditions as follows: we assume that there is no flow of energy through the end of the line where $z=1$, or in other words $i(1,t)=0$. In addition assume that $E=0$ in (\ref{eq_b1}) (the boundary conditions at $z=0$) This results in the boundary condition $v_0 = i_0 R_0$ at $z=0$. We make use of the new ``admissible pairs" $(\tilde A, \tilde P)$ for the transmission line derived as in equations(\ref{eq_newP},\ref{eq_newA}) resulting in the mixed potential function
\begin{align}
\mathcal P = \tilde P + P^0
\end{align}
where $\tilde P$ is defined as in (\ref{eq_newP}) and $P^0 = \frac{1}{2}i_0^2 R_0$, taking care of the boundary resistance at $z=0$. Differentiating $\mathcal P$ along the trajectories of the transmission line results  in
\begin{align}
\dot {\mathcal P} \le v_0 \frac{di_0}{dt} = u_0 y_0
\end{align}
resulting in passivity with respect to the voltage and the derivative of the current at $z=0$.
In deriving this new passivity property we have made use of the fact that a purely resistive network defines a passive system with respect to $(i_0,\frac{dv_0}{dt})$ as well as with respect to $(v_0, \frac{di_0}{dt})$.

\section{Control by Interconnection}
In this section we consider stabilization of the system at a nontrivial equilibrium point via boundary control. The argument used here is same
as that presented in \cite{hugo}, where the authors have presented a boundary control law for a mixed finite and infinite dimensional system via energy shaping methods.
In this section we propose a boundary control law via shaping the power of the infinite-dimensional system. As in the finite dimensional case, the method relies on finding Casimir functions for the closed-loop system.
Consider the equations of the transmission line with a new admissible pair, written as
\begin{align}
\frac{\partial }{\partial t}
\begin{bmatrix}
i(z,t) \\
v(z,t)
\end{bmatrix}
=
\begin{bmatrix}
-\frac{1}{L} & -\frac{2}{GL} \frac{\partial}{\partial z} \\
0 & -\frac{1}{C}
\end{bmatrix}
\begin{bmatrix}
\delta_i { \mathcal P} \\
\delta_v { \mathcal P}
\end{bmatrix}
\label{eq_newad}
\end{align}
with ${\mathcal P}$ defined as in (\ref{eq_newP}). For simplicity, let us assume that there is no power flowing through the end of the transmission line at $z=1$ (i.e $(i(1,t) = 0$). In this case
\[
\frac{d}{dt}  {\mathcal P} \le v_0 \frac{di_0}{dt} = u_0 y_0
\]
At the end of the transmission line where $z=0$, interconnect the system with a controller of the form
\begin{align}
\dot \xi = u_c, ~~ y_c = \frac{\partial H_c(\xi)}{\partial \xi}
\label{eq_ctrl}
\end{align}
where $\xi, u_c$ and $y_c$ are respectively the state, input and output of the controller. $H_c(\xi)$ denotes the power function of the controller.
\begin{align}
\begin{bmatrix}
u_0  \\
u_c
\end{bmatrix}
=
\begin{bmatrix}
0 & -1 \\
1 & 0
\end{bmatrix}
\begin{bmatrix}
y_0 \\
y_c
\end{bmatrix}
\label{eq_int}
\end{align}
This results in the following ``closed-loop" system
\begin{align}
\begin{bmatrix}
\dot \xi \\
\frac{di(z,t)}{dt} \\
\frac{d v(z,t)}{dt}
\end{bmatrix}
=
\begin{bmatrix}
0 & -\frac{1}{L }~\cdot \mid_{z=0} & -\frac{2}{GL} \frac{\partial }{\partial z} ~ \cdot \mid_{z=0} \\
0  & -\frac{1}{L} & -\frac{2}{GL} \frac{\partial }{\partial z} \\
0 & 0 & -\frac{1}{C}
\end{bmatrix}
\begin{bmatrix}
\frac{\partial H_c}{\partial \xi} \\
\delta_i  {\mathcal P} \\
\delta_v  {\mathcal P}
\end{bmatrix}
\label{eq_cl}
\end{align}
where $\cdot \mid_{z=0}$ denoted the value of the the corresponding variational derivative of $ {\mathcal P}$ evaluated at $z=0$. In the above expression the first equation would mean
\[
\dot \xi = -\frac{1}{L} \delta_i  {\mathcal P}\mid_{z=0} - \frac{2}{GL}\frac{\partial}{\partial z} \left (\delta_v  {\mathcal P}\right )\mid_{z=0}
\]
To achieve the control objectives we look for Casimir functions for the closed-loop system in (\ref{eq_cl}). A functional $\mathcal C (\xi, i,v)$ is a Casimir functional provided
\begin{align}
\left (
\begin{bmatrix}
0 & -\frac{1}{L }~\cdot \mid_{z=0} & -\frac{2}{GL} \frac{\partial }{\partial z} ~ \cdot \mid_{z=0} \\
0  & -\frac{1}{L} & -\frac{2}{GL} \frac{\partial }{\partial z} \\
0 & 0 & -\frac{1}{C}
\end{bmatrix}
\begin{bmatrix}
\frac{\partial \mathcal C}{\partial \xi} \\
\delta_i  {\mathcal C} \\
\delta_v  {\mathcal C}
\end{bmatrix} \right )^T = 0
\end{align}
From the above expression we can see that $\delta_v \mathcal C = 0$, which implies that the Casimir functions are independent of $v$. In addition we can also infer that
$\frac{\partial}{\partial z} \delta_i \mathcal C = 0$, which means that $ \delta_i \mathcal C$ does not vary with the spatial parameter $z$. Further if we restrict the Casimirs to functionals of the form $\mathcal C = -\xi + \mathcal F (i(z,t))$, then the conditions on the functional $\mathcal F(i(z,t))$ reduce to $\delta_i \mathcal F = 1$. While deriving the above conditions we have used the fact that the conjugate of the $\frac{\partial}{\partial z}$  (skew-symmetric) operator is -$\frac{\partial}{\partial z}$.

\subsection{Control Design and Stability Analysis }
The control objective here is to stabilize the the voltages and currents of the transmission line to steady state values $(v^{\ast}, i^{\ast})$ such that
\[
Gv^{\ast} + \frac{\partial i}{\partial z} \mid_{i=i^{\ast}} = 0.
\]
From the above discussions it follows that any functional of the form $-\xi + \int_0^1 i(z,t)dz =  \text{constant} $ is a Casimir for the closed-loop system. Then in the invariant set
\[
\Omega = \{ \mathcal X \mid \xi = \int_0^1 i(z,t)dz + c \}
\]
the closed loop mixed potential function equals
\[
\mathcal P_d = {\mathcal P} + H_c({\xi}) + c, ~ \text{with } c \text{ a constant}
\]
We need to choose a $H_c$ in such a way that the closed loop mixed potential function $\mathcal P_d $ satisfies the conditions (a)-(c) of Theorem \ref{thm}.
We do this as follows: Chose $H_c(\xi)$ as
\begin{align}
H_c(\xi) =
\int_Z\left(\dfrac{K}{2}(i-i^{\ast})^2-ii^{\ast}R+\frac{\partial{i}}{\partial{x}}v^{\ast}\right)dz
\end{align}
The resulting  $\mathcal P_d $ is
\begin{align}
 & \mathcal P_d  =   \int_0^1\left ( \frac{1}{2} Ri^2 + \frac{1}{2} Gv^2 + v \frac{\partial i}{ \partial z} + \frac{1}{G} \left ( \frac{\partial i}{\partial z}\right )^2\right )dz \nonumber \\
 & +
 \int_0^1 \left(\dfrac{K}{2}(i-i^{\ast})^2-ii^{\ast}R+\frac{\partial{i}}{\partial{z}}v^{\ast}\right)dz
 + c \nonumber \\
 & =\int_0^1
\left(\dfrac{1}{2G}\left(Gv+\frac{\partial{i}}{\partial{z}}\right)^2+\dfrac{1}{2G}\left(\frac{\partial{i}}{\partial{z}}\right)^2\right)dz\nonumber \\
&+\int_0^1 \left(\dfrac{K}{2}(i-i^{\ast})^2+\frac{R}{2}\left(i-i^{\ast}\right)^2+\frac{\partial{i}}{\partial{z}}v^{\ast}-\frac{R}{2}i^{\ast2}\right)dz \nonumber \\
& +c\nonumber \\
& =\dfrac{1}{2G}\int_0^1
\left(\left(Gv+\frac{\partial{i}}{\partial{z}}\right)^2+G(K+R)(i-i^{\ast})^2+\right. \nonumber\\
\nonumber
&\left.\left(Gv^{\ast}+\frac{\partial{i}}{\partial{z}}\right)^2\right)dz-\dfrac{1}{2}\int_Z\left(Ri^{\ast2}+Gv^{\ast2}\right)dz+c\nonumber\\
%
%=  &  \int_0^1 \left (\frac{1}{2G} \left ( Gv + \frac{\partial i}{\partial z}\right )^2 + \frac{K+R}{2} (i-i^{\ast})^2 \right . \nonumber \\
%& \left . \frac{1}{2G}\left ( \frac {\partial i}{\partial z} \right )^2 - i\frac{\partial v}{\partial z} \mid_{v=v^{\ast}} \right ) - \frac{R}{2} \int_0^1 {i^{\ast}}^2 dz + c
%
\end{align}
 Finally choosing $c =   \frac{1}{2} \int_Z (R{i^{\ast}}^2 + G {v^{\ast}}^2)dz $, we arrive at
\begin{align}
& \mathcal P_d= \dfrac{1}{2G}\int_Z
\left(\left(Gv+\frac{\partial{i}}{\partial{z}}\right)^2+G(K+R)(i-i^{\ast})^2+\right.\\
\nonumber
&\left.\left(Gv^{\ast}+\frac{\partial{i}}{\partial{z}}\right)^2\right)dz\nonumber
\end{align}
Such a $\mathcal P_d$ satisfies the conditions 3a, 3b of Theorem \ref{thm}.
The closed loop dynamics with this new $\mathcal P_d$, the ``shaped power" takes the form
\begin{align}
\frac{\partial }{\partial t}
\begin{bmatrix}
i(z,t) \\
v(z,t)
\end{bmatrix}
=
\begin{bmatrix}
-\frac{1}{L} &- \frac{2}{GL} \frac{\partial}{\partial z} \\
0 & -\frac{1}{C}
\end{bmatrix}
\begin{bmatrix}
\delta_i { \mathcal P_d} \\
\delta_v { \mathcal P_d}
\end{bmatrix}
\end{align}
Furthermore it can be verified that, $ \delta_u \tilde {\mathcal P_d}= 0$ at the desired equilibrium, which proves the condition for the desired equilibrium.
Lastly, under the condition that $\sqrt{C/L}/G ||\frac{\partial}{\partial z}|| \le 1$ we can also show that
\[
\dot{ \tilde{\mathcal P_d}} \le 0
\]
satisfying condition 3c of Theorem \ref{thm}.
We have thus proved the following:
\begin{proposition}
Consider the transmission line whose dynamics are given by (\ref{eq_newad}), together with the controller of the form (\ref{eq_ctrl}), interconnected to the plant (the transmission line), via the constraints (\ref{eq_int}), and the controller power function given by (22). The resulting closed-loop system has a stable equilibrium at $(v^{\ast}, i^{\ast})$ where
$
Gv^{\ast} + \frac{\partial i}{\partial z} \mid_{i=i^{\ast}} = 0.
$
\end{proposition}

\section{remarks on other alternative passive maps}
The earlier work on passivity based control of infinite-dimensional systems relied on using Hamiltonian as the Lyapunov functional and the passivity property was proven by means of boundary port variables which are power conjugate. In case of transmission line the boundary variables, which prove the passivity property are the voltages and currents at the boundary. In the previous section we have shown how by the use of the mixed-potential function as the storage functional, resulting in a new passivity property with the voltages and the derivatives as the current at the boundary as the new input-output port variables. The storage functional in this case has the dimensions of power, that is energy/time. Similarly there exist possibilities of deriving new passivity properties by using a storage functional which has dimensions of power/time and the derivatives of both the voltages and current as the new port variables. This can be shown as follows:
%\section{distributed network}
%%
%\begin{itemize}
%%
%\item $L\frac{\partial{i}}{\partial{t}}$ =
%$-\frac{\partial{v}}{\partial{x}}-Ri$
%\item $C\frac{\partial{v}}{\partial{t}}$ =
%$-\frac{\partial{i}}{\partial{x}}-Gv$
%\item $x=0$ $\Rightarrow$   \quad $E=v_0$
%
%\item $x=1$ $\Rightarrow$  \quad $I_s=i_1$
%\end{itemize}
Consider a functional given by  $P=\lambda P_1+P_2$ (also refer \cite{BM03}), where
\beqn
P_1&=& \int_0^1
\left(-\frac{1}{2}Ri^2+\frac{1}{2}Gv^2-i\frac{\partial{v}}{\partial{z}}\right)dz,
\eeqn

\beqn
P_2&=&
\int_0^1L\left(\frac{\partial{i}}{\partial{t}}\right)^2+C\left(\frac{\partial{v}}{\partial{t}}\right)^2dz.
\eeqn
Their differentials with time are given as:
\beq
\dot{P}_1&=&L\int_0^1\left(\frac{\partial{i}}{\partial{t}}\right)^2dz-C\int_0^1\left(\frac{\partial{v}}{\partial{t}}\right)^2dz
%+\frac{\partial{\left(v_1.i_1\right)}}{\partial{t}}
\label{pdot}
\eeq
\beq
\dot{P}_2&=&\int_0^1L\left(\frac{\partial{i}}{\partial{t}}\right)\frac{\partial{}}{\partial{t}}\frac{\partial{i}}{\partial{t}}+C\left(\frac{\partial{v}}{\partial{t}}\right)\frac{\partial{}}{\partial{t}}\frac{\partial{v}}{\partial{t}}dz\nonumber\\
&=&-\int_0^1R\left(\frac{\partial{i}}{\partial{t}}\right)^2dz-\int_0^1G\left(\frac{\partial{v}}{\partial{t}}\right)^2dz \nonumber \\
&& -\left. \left(\frac{\partial{v}}{\partial{t}}.\frac{\partial{i}}{\partial{t}}\right)\right|_0^1\nonumber\\\label{qdot}
\eeq

Finally from \label{pdot}, \label{qdot} We get
\beq
\dot{\hat{P}}&=&\lambda\dot{P}_1+\dot{P}_2\nonumber\\
&=&-(R-\lambda
L)\int_0^1\left(\frac{\partial{i}}{\partial{t}}\right)^2dz\nonumber\\&&-(\lambda
C+G)\int_0^1\left(\frac{\partial{v}}{\partial{t}}\right)^2dz-\left.\left(\frac{\partial{v}}{\partial{t}}.\frac{\partial{i}}{\partial{t}}\right)\right|_0^1\nonumber\\
& \le&  \left.\left(\frac{\partial{v}}{\partial{t}}.\frac{\partial{i}}{\partial{t}}\right)\right|_0^1 \label{qdot}
\eeq
provided $0<\lambda < \frac{R}{L}$\\
Note that $P_1$ by itself has the dimensions of power whereas $P_2$ has dimensions of power/time. We thus multiply $P_1$ with a constant $\lambda$ which has units of 1/time such that the new storage functional $P = P_1 + P_2$ has dimensions of power/time and preserves the passivity property by using the port variables as the derivatives of both the voltages and currents at the boundary.
\addtolength{\textheight}{-16.3cm} 
\section{Conclusions and Future work}
In this paper we present an example showcasing power shaping for infinite-dimensional systems within the BM framework. We also present case for generation of other alternate passive maps. Future work would focus on generalizing these concepts to a general class of infinite-dimensional systems from various physical domains. Also it would be of interest to analyze the underlying geometric structure of  infinite-dimensional systems under the BM framework.

 \section*{Acknowledgements}
 The first author would like to thank Pujita Raman for her contribution during the initial phase of this work.


\begin{thebibliography}{00}


\bibitem{BM01}R. K. Brayton and J. K. Moser. A theory of nonlinear networks I. {\em Q.
Appl. Math.}, vol. 22, no. 1, pp. 1-33, Apr. 1964.

\bibitem{BM02}R. K. Brayton and J. K. Moser. A theory of nonlinear networks II. {\em Q.
Appl. Math.}, vol. 22, no. 2, pp. 81-104, Jul. 1964.

\bibitem{BM03} R. K. Brayton and W. L. Miranker. A stability theory for nonlinear
mixed-initial boundary value problems, {\em Arch. Rat. Mech. Anal.}, vol. 17,
no. 5, pp. 358Ð376, Dec. 1964.

\bibitem{geoplex} V. Duindam, A. Macchelli, S. Stramigioli and H. Bryuninckx, Eds. \newblock {\it Modeling and Control of Complex Physical Systems:
The port-Hamiltonian approach.} \newblock Springer, 2009.

\bibitem{ifac} E. Garcia- Canseco, R. Pasumarthy, A.J. van der Schaft and R. Ortega.
\newblock On control by interconnection of port-Hamiltonian systems. {\em Proc. IFAC world Congress}, July 2005.

\bibitem{guido} Guido Blankenstein.
\newblock Power Balancing For A New Class of Non-Linear Systems and Stabilization of RLC circuits. {\em International Journal of Control},
Vol. 78, Issue 3, pp. 159 -- 171, February 2005.


\bibitem{eloisa} Eloisa Garcia-Canseco and Romeo Ortega.
\newblock A new passivity property of linear RLC circuits with application to Power Shaping
Stabilization. {Proc. American Control Conference}, 2004.

\bibitem{ejos} E. Garcia Canseco, D. Jeltsema, J.M.A. Scherpen, and R. Ortega. Power-based control of physical systems, {\em Automatica}, 46, pp. 127-132, 2010.

\bibitem{jos} D. Jeltsema, R. Ortega and J.M.A. Scherpen. \newblock On Passivity and Power-Balance Inequalities of Nonlinear RLC Circuits, {\em IEEE Transactions on Circuits and Systems-I: Fundamental Theory and Applications}, Vol. 50, No. 9, 2003, pp. 1174-1179.

\bibitem{josh}D. Jeltsema, R. Ortega and J.M.A. Scherpen. \newblock An Energy-Balancing Perspective of Interconnection and Damping Assignment Control of Nonlinear Systems, {\em Automatica}, Vol. 40, No. 9, September 2004, pp. 1643-1646.

\bibitem{ja}D. Jeltsema and A. J. van der Schaft. Pseudo-gradient and Lagrangian
boundary control formulation of electromagnetic fields,  {\em J. Phys. A:
Math. Theor.}, vol. 40, pp. 11627Ð11643, 2007.

\bibitem{jeltcsm} D. Jeltsema and J.M.A. Scherpen. Multidomain
\newblock Multidomain Modeling of Nonlinear Networks and Systems Energy- and Power-based perspectives.
{\em IEEE Control Systems Magazine}, pp 28 -- 59, August 2009.

\bibitem{lzm} Y. Le Gorrec, H. Zwart and B.M.J. Maschke. \newblock Dirac structures and boundary control systems
associated with skew-symmetric differential operators. 
\newblock {\em SIAM Journal on Control
and Optimization}, 44 (5). pp. 1864-1892, 2005.

\bibitem{ortegapower} R. Ortega, D. Jeltsema, and J. M. A. Scherpen.
\newblock Power Shaping: A New Paradigm for Stabilization of Nonlinear RLC Circuits.
{\em IEEE Transactions on Automatic Control}, VOL. 48, NO. 10, pp. 1762 -- 1767, OCTOBER 2003

\bibitem{ortaut} R. Ortega, A. J.  van der Schaft, B.M. Maschke and G Escobar.
\newblock Interconnection and damping assignment
passivity-based control of port-controlled Hamiltonian
systems. {\em Automatica}, Vol. 38, Issue 4,  pp. 585-596, 2002.

\bibitem {ramijc} R. Pasumarthy and A.J. van der Schaft. \newblock Achievable Casimirs and its implications on control
by interconnection of port-Hamiltonian systems. \newblock {\it International Journal of Control,} vol. 80, no.9, pp. 142--1438, 2007.

\bibitem{rammtns} R. Pasumarthy and A.J. van der Schaft. \newblock On interconnections of infinite-dimensional port-Hamiltonian systems. 
\newblock {\em Proceedings 16th International Symposium on Mathematical Theory of Networks and Systems (MTNS)}, July 2004, Leuven, Belgium.

\bibitem{hugo} H. Rodriguez, A.J. van der Schaft and R. Ortega. On
stabilization of nonlinear distributed
parameter port-controlled Hamiltonian systems via energy shaping. In {\it
Proceedings
of the 40th IEEE conference on decision and control, Orlando}, FL, December
2001.

\bibitem {l2gain}A.J. van~der Schaft.
\newblock {\em $L_2$-{G}ain and {P}assivity {T}echniques in {N}onlinear
{C}ontrol}. \newblock Springer-Verlag, 2000.

\bibitem {stokes}A.J. van~der Schaft and B.M. Maschke. \newblock Hamiltonian
formulation of distributed-parameter systems with boundary energy flow.
\newblock {\it Journal of Geometry and
Physics, vol.42}, pp.166-194, 2002.



 \end{thebibliography}
\end{document}